\documentclass[a4paper,12pt,leqno]{article}
\usepackage{latexsym}
\usepackage[all]{xy}

\usepackage{amssymb} 
\usepackage{amsmath} 
\usepackage{theorem}

\def\Z{{\mathbb{Z}}}
\def\Q{{\mathbb{Q}}}
\def\K{{\mathbb{K}}}
\def\CC{{\mathbb{C}}}
\def\R{{\mathbb{R}}}
\def\CC{{\mathbb{C}}}
\def\A{{\mathcal{A}}}
\def\B{{\mathcal{B}}}

\DeclareMathOperator{\codim}{codim}

\DeclareMathOperator{\Der}{Der}

\DeclareMathOperator{\soc}{soc}
\DeclareMathOperator{\socdeg}{socdeg}

\numberwithin{equation}{section}

\newcommand{\owari}{\hfill$\square$}

\theoremstyle{break}
\newtheorem{theorem}{Theorem}[section]
\newtheorem{prop}[theorem]{Proposition}
\newtheorem{cor}[theorem]{Corollary}
\newtheorem{lemma}[theorem]{Lemma}
\newtheorem{define}[theorem]{Definition}
\newtheorem{rem}[theorem]{Remark}

\newtheorem{example}[theorem]{Example}
\newtheorem{problem}[theorem]{Problem}
\newtheorem{conj}[theorem]{Conjecture}
\newtheorem{question}[theorem]{Question}

\title{Solomon-Terao algebra of hyperplane arrangements}
\author{
Takuro Abe, 
Toshiaki Maeno, Satoshi Murai and Yasuhide Numata
}

\date{\today} 

\begin{document}

\maketitle

\begin{abstract}
We introduce a new algebra associated with a hyperplane arrangement $\A$, called the 
Solomon-Terao algebra $ST(\A,\eta)$, where $\eta$ is a 
homogeneous polynomial. It is shown by Solomon and Terao that 
$ST(\A,\eta)$ is Artinian when $\eta$ is generic. 
This algebra can be considered as a generalization of coinvariant algebras in the 
setting of hyperplane arrangements. The class of Solomon-Terao algebras 
contains cohomology rings of regular nilpotent Hessenberg varieties.
We show that $ST(\A,\eta)$ is a complete intersection if and only if 
$\A$ 
is free. 
We also give a factorization formula of the Hilbert polynomials when $\A$ is free, and 
pose several related questions, problems and conjectures.
\end{abstract}

\section{Introduction}

The aim of this article is to introduce a new algebra, called the 
\textbf{Solomon-Terao 
algebra} and the \textbf{Solomon-Terao complex} associated with hyperplane arrangements. The classical and well-studied algerbra 
of hyperplane arrangement is the logarithmic derivation module, and our Solomon-Terao 
algebra is defined by using logarithmic derivation modules. 
The Solomon-Terao algebra has two remarkable aspects. The first one is that
it corresponds to a specialization of the Solomon-Terao polynomial defined in \cite{ST},
while the famous Orlik-Solomon algebra in \cite{OS} reflects another kind of specialization of
the Solomon-Terao polynomial. Hence the Solomon-Terao algebra is considered to be
comparable with the Orlik-Solomon algebra, which is isomorphic to 
the cohomology ring of the complement of the hyperplane arrangement 
when the base field is $\CC$. 
We study the algebraic structure of the Solomon-Terao 
algebra in \S1.1. 

The second aspect is a geometric feature of the Solomon-Terao algebra,
which gives some support of the suitability of our definition.
The Solomon-Terao algebra happens to be isomorphic to the cohomology ring of
some varieties, analogously to the Orlik-Solomon algebra.
In a typical case, the Solomon-Terao algebra can be isomorphic to the cohomology ring
of the flag variety or the coinvariant algebra of the reflection group. 
This part is described in \S1.2.


\subsection{Solomon-Terao algebra and main results}

Let us introduce 
several definitions.
Let $\K$ be an algebraically closed field, $V=\K^\ell$ and 
$S:=\mbox{Sym}(V^*)$ its coordinate ring. Let us fix a coordinate system $x_1,\ldots, x_\ell$ for 
$V^*$ such that $S=\K[x_1,\ldots,x_\ell]$. The $\K$-linear $S$-derivation module $\Der S$ is 
a rank $\ell$ free module defined by 
$$
\Der S:=\oplus_{i=1}^\ell S \partial_{x_i}.
$$
Also, let $\Der^p S:=\wedge^p \Der S\ (p \ge 0)$, agreeing that $\Der^0 S=S$. 
Let $\A$ be an arrangement of linear hyperplanes, i.e., a finite set of linear hyperplanes 
in $V$. For each $H \in \A$ fix a linear form $\alpha_H \in V^*$ such that $\ker \alpha_H=H$. Let $Q(\A):=\prod_{H \in \A} \alpha_H$. 
Now we can define the logarithmic derivation modules $D^p(\A)$ for $\A$ as follows:

\begin{define}
For $p \ge 0$, define 
$$
D^p(\A):=\{\theta \in \Der^p S \mid 
\theta(\alpha_H,f_2,\ldots,f_p) \in S \alpha_H\ (\forall H \in \A,\ 
\forall f_2,\ldots,f_p \in S)\}.
$$
\label{log}
\end{define}

The logarithmic derivation module was introduced by K. Saito for the study of the universal unfolding of the isolated hypersurface singularity, see 
\cite{S} for example. The logarithmic derivation module has been studied mainly for $p=1$
particularly in case of hyperplane arrangements. 
On the other hand, by using $D^p(\A)$ for all $p$, Solomon and Terao defined the following 
interesting series.

\begin{define}[\cite{ST}, Section One]
Define the \textbf{Solomon-Terao polynomial} $\Psi(\A;x,t)$ by 
$$
\Psi(\A;x,t):=t^\ell\sum_{p=0}^\ell 
\mbox{Hilb}(D^p(\A);x)(\frac{1-x}{t}-1)^p.
$$
\label{STpolynomial}
\end{define}

Here for the graded $S$-module $M$, 
$$
\mbox{Hilb}(M;x):=\sum_{i=0}^\infty (\dim_\K M_i) x^i
$$
is the Hilbert series of $M$. 
In the definition above, the Solomon-Terao polynomial
seems to be introduced as a series. However, in fact 
it is a polynomial.

\begin{theorem}[\cite{ST}, Proposition 5.3]
$\Psi(\A;x,t) \in \Q[x,t]$.
\label{polynomial}
\end{theorem}

Moreover, we have the following astonishing result by \cite{ST}.

\begin{theorem}[\cite{ST}, Theorem 1.2]
Let $\K=\CC$ and $\pi(\A;t)$ the topologocal Poincar\`{e} 
polynomial of $M(\A):=V \setminus \cup_{H \in \A} H$. Then 
$$
\Psi(\A;1,t)= \pi(\A;t).
$$
\label{STmain}
\end{theorem}

In fact, $\pi(\A;t)$ can be defined over an arbitrary field $\K$ by using 
combinatorial data of $\A$ see \S2. 
Hence Theorem \ref{STmain} connects algebra, topology and 
combinatorics of $\A$. By \cite{OS}, we know that there is the algebra $A(\A)$ called the 
\textbf{Orlik-Solomon algebra} depending only on the intersection lattice $L(\A)$ such that 
$$
A(\A) \simeq H^*(M(\A) ,\Z)
$$
when $\K=\CC$, see \cite{OT} \S3 and \S5 for details. In particular, Theorem \ref{STmain} implies that 
\begin{equation}
\Psi(\A;1,t)=\pi(\A;t)=\mbox{Hilb}(A(\A) \otimes \Q;x).
\label{OS}
\end{equation}

As we see that the specialization $\Psi(A;1,t)$ has
a geometric meaning in Theorem \ref{STmain}, it is natural to ask whether the specialization
with respect to the $t$-variable has a nice interpretation. 
For example, can we understand $\Psi(\A;x,1)$ 
by using algebra, geometry of arrangements or other geometric objects? In this subsection 
we give an answer to this problem from algebraic point of view. 
Let us introduce algebraic 
counter part of $\Psi(\A;x,1)$ in the following.

\begin{theorem}[\cite{ST}]

Let $d$ be a non-negative integer 
and $S_d$ the set of all homogeneous polynomials of degree $d$. Fix 
$\eta \in S_d$.
Also, define the boundary map $\partial_p:D^p(\A) \rightarrow 
D^{p-1}(\A)\ (p=1,\ldots,\ell)$ by 
$$
\partial_p(\theta)(f_2,\ldots,f_p):=
\theta(\eta,f_2,\ldots,f_p)
$$
for all $f_2,\ldots,f_p \in S$. 
We call the complex $(D^*(\A),\partial_*)$ the \textbf{Solomon-Terao complex} 
of degree $d$ with respect to $\eta \in S_d$. Define their cohomology group
$$
H^p(D^*(\A),\partial_*):=\mbox{ker} \partial_p / \mbox{Im} \partial_{p+1}.
$$
Then 

(1)\,\, 
there is a non-empty Zariski open set $U_d =U_d(\A) \subset S_d$ such that 
every cohomology of the Solomon-Terao complex with respect to 
$\eta \in U_d$ is of finite dimension over $\K$. For the details of $U_d$, see \S 2. 

(2)\,\,
If $\mbox{pd}_S D^p(\A) \le \ell-p$ for all $p=1,2,\ldots,\ell$ (such an arrangement is called 
\textbf{tame}), then 
$H^i(D^*(\A),\partial_*)=0$ for $i \neq 0$.
\label{ST}
\end{theorem}

%

\begin{define}[Solomon-Terao algebra]
In the notation of Theorem \ref{ST}, define 
$ST(\A,\eta):=H^0(D^*(\A),\partial_*)$ and let us call $ST(\A,\eta)$ 
the \textbf{Solomon-Terao algebra} of $\A$ with respect to $\eta$. 
We call $\mathfrak{a}(\A,\eta):=\{ \theta (\eta) \in S \mid \theta \in D(\A)\}=
\mbox{Im} \partial_1$ the 
\textbf{Solomon-Terao ideal} of $\A$ with respect to $\eta$, i.e., $S/\mathfrak{a}(\A,\eta)=ST(\A,\eta)$. 
\label{STalgebra}
\end{define}

\begin{rem}
By definition, the structure of the Solomon-Terao algebra depends on the choice of the 
polynomial $\eta \in U_d(\A)$. See Example \ref{notisom} for details.
\end{rem}

The Solomon-Terao algebra can be defined for all arrangements, but the most 
useful case is when $\A$ is tame. In fact, we can show that 
the Solomon-Terao algebra is the algebraic counterpart of 
$\Psi(\A;x,1)$ when $\A$ is tame.

\begin{theorem}
Let $\A$ be tame and $\eta \in U_2(\A)$. Then we have  
$$
\mbox{Hilb}(ST(\A, \eta);x)=\Psi(\A;x,1).
$$
In particular, 
$\mbox{Hilb}(ST(\A, \eta);1)=\pi(\A;1)$ coincides with the number of chambers when 
$\K=\R$, and with the total Betti numbers of $M(\A)$ when $\K=\CC$. 
\label{poin2}
\end{theorem}

Theorem \ref{poin2} is essentially proved in \cite{ST}. 
We have reformulated the result in form of Theorem \ref{poin2} to explain a reason
to consider the Solomon-Terao algebra. 
Theorem \ref{poin2} affords a good motivation to study. 
However, that is 
the motivation to define and study $ST(\A,\eta)$. 
Hence the Solomon-Terao algebra is closely related to the Solomon-Terao 
polynomial from algebraic point of view when it is tame.
Though tameness is a generic property, to check the tameness for a given arrangement is very hard. 
Fortunately, one of the most famous classes of hyperplane arrangements is known to be tame.

\begin{define}
An arrangement $\A$ is \textbf{free} with $\exp(\A)=(d_1,\ldots,d_\ell)$ if 
$D(\A)$ is a free $S$-module of rank $\ell$ with homogeneous basis $\theta_1,\ldots,\theta_\ell,\ 
\deg \theta_i=d_i\ (i=1,\ldots,\ell)$.
\label{free}
\end{define}

When $\A$ is free, $D^p(\A)$ is also free (see \cite{OT} or \cite{ST}). Thus the 
freeness implies the tameness. 
Since the freeness is a very strong property of arrangements, it is worth studying 
$ST(\A,\eta)$ when $\A$ is free, which is our second main result. 
To state it, let us recall some fundamental definitions on 
commutative ring theory. Let $M$ be a graded $S$-module. Let 
$M_n$ denote the homogeneous degree $n$-part of $M$. Then 
the \textbf{socle} $\soc(M)$ of $M$ is defined as 
$$
\soc(M):=0:_M S_+,
$$
where $S_+=(x_1,\ldots,x_\ell)$. When $S/I$ is a Gorenstein 
$\K$-algebra for an ideal $I \subset S$, 
$\dim_\K \soc (S/I)=1$, hence there is an integer $r$ such that 
$\soc(S/I)=(S/I)_r$. 
We call $r$ the \textbf{socle degree} of $S/I$ and denote
it by $\socdeg(S/I)$.
%
Now we can state the second 
main theorem in this article.

\begin{theorem}[Freeness and C.I.]
(1)\,\,
Assume that $\A$ is free with $\exp(\A)=(d_1,\ldots,d_\ell)$ and let $\eta \in U_d$. Then 
$ST(\A, \eta)$ is a complete intersection with
$$
\mbox{Hilb}(ST(\A,\eta);x)=\prod_{i=1}^\ell (1+x+\cdots+x^{d_i+d-2}).
$$
Hence $\socdeg ST(\A,\eta)=|\A|+\ell(d-2)$, 
and 

(2)\,\,
conversely, 
if $ST(\A,\eta)$ is a complete intersection, then $\A$ is free. If 
$$
\mbox{Hilb}(ST(\A,\eta);x)=\prod_{i=1}^\ell (1+x+\cdots+x^{e_i}),
$$
then $\exp(\A)=(e_1-d+2,\ldots,e_\ell-d+2)$.

\label{factorization}
\end{theorem}

\begin{rem}
When $d=2$, Theorem \ref{factorization} is also proved by 
Epure and Schulze in \cite{ES} independently. Explicitly, when $d=2$, they 
proved the same result 
not only for hyperplane arrangements but also for hypersurface singularities. 
\end{rem}

Theorems \ref{poin2} and \ref{factorization} helps us to investigate the algebraic 
structure of $ST(\A,\eta)$ in terms of commutative ring theory. Then the next question 
is whether we have a nice geometric understanding of $ST(\A,\eta)$ and 
$\Psi(\A;x,1)=\mbox{Hilb}(ST(\A,\eta);x)$ when $\A$ is tame. Let us give an answer 
from classical result by Borel in \cite{B}, and the recent results on Hessenberg 
varieties in \cite{AHMMS} in the next subsection.


\subsection{Geometry of the Solomon-Terao algebra}

In this subsection, we show the relation between $ST(\A,\eta)$ and 
the cohomology ring of some variety, which is analogous to the one 
between the Orlik-Solomon algebra $A(\A)$ and 
the open manifold $M(\A)$ when $\K=\CC$. In this subsection let $\K=\CC$. 

First let $W$ be the irreducible crystallographic Weyl group acting on $V$. 
Let $G$ be the corresponding complex semisimple linear algebraic group, and 
$B$ the fixed Borel subgroup.  
Let 
$\A=\A_W$ be a set of reflecting  
hyperplanes of all reflections of the Weyl group $W$ (so called the \textbf{Weyl arrangement}). 
By the result of K. Saito (see \cite{S} for example), $\A_W$ is free with exponents 
$(d_1^W,\ldots,d_\ell^W)$ coinciding with the exponents of $W$. Combining the results in \cite{B} and \cite{ST}, we have 
\begin{equation}
\Psi(\A_W;x,1)=\mbox{Poin}(G/B;\sqrt{x})=\prod_{i=1}^\ell (1+x+\cdots+x^{d_i^W}).
\label{Weyl}
\end{equation}
Here $G/B$ is the flag variety corresponding to $W$. 
Hence $\Psi(\A_W;x,t)$ has two important specializations for $x=1$ and $t=1$ in the geometric 
point of view. Also, let $S^W$ denote the $W$-invaraint part of the polynomial ring with the $W$-action, and let 
$\mbox{coinv(W)}:=S/(S^W_+)$ the coinvariant algebra. Then it is well-known that 
\begin{equation}
\mbox{Hilb}(\mbox{coinv}(W);x)=\mbox{Poin}(G/B;\sqrt{x})=\Psi(\A_W;x,1).
\label{coinvHinb}
\end{equation}
Hence when $\A=\A_W$, the algebraic counterpart of $t=1$ is the coinvariant algebra, which is also known to be isomorphic to the cohomology ring of the flag variety $G/B$ by Borel in 
\cite{B}. In fact, 
we get a natural interpretation 
of the Solomon-Terao algebra for $\A=\A_W$. 
Let 
$P_1$ be the lowest degree basic invariant of $S^W$. Then Theorem 3.9 in \cite{AHMMS} 
shows that 
$$
ST(\A_W,P_1) \simeq \mbox{coinv}(W) \simeq H^*(G/B,\CC).
$$
Thus we can understand the Solomon-Terao algebra from geometric point of view 
in a way suggestive of 
the Orlik-Solomon algebra when $\A=\A_W$. The above isomorphism is now extended to a wider class. We refer to the results in \cite{AHMMS}.

%

\begin{define}
Let $\Phi$ be the root system with respect to $W$ and fix a 
positive system $\Phi^+$. 
Let $I \subset \Phi^+$ be a \textbf{lower ideal}, i.e., the set satisfying 
that, if $\beta \in I,\ \gamma \in \Phi^+$ and 
$\beta-\gamma \in \sum_{i=1}^\ell \Z_{\ge 0} \alpha_i$ for 
the simple system $\alpha_1,\ldots,\alpha_\ell$ of $\Phi^+$, then $\gamma \in I$. Then 
$\A_I:=\{\alpha=0\mid \alpha \in I\}$ is called the \textbf{ideal arrangement}. 
\label{ideal}
\end{define}

The freeness of the ideal arrangements are known as follows:

\begin{theorem}[Theorem 1.1, \cite{ABCHT}]
Let $I \subset \Phi^+$ be a lower ideal. Then $\A_I$ is free with exponents $
(d_1^I,\ldots,d_\ell^I)$ which coincides with the dual partition of the height distribution of 
the positive roots in $I$. 
\label{idealarr}
\end{theorem}

Hence Theorem \ref{factorization} is applicable to the algebra $ST(\A_I,P_1)$. 
On the other hand, we can also associate a variety with the lower ideal, so called the 
regular nilpotent Hessenberg variety, see \cite{DPS} or 
\cite{AHMMS} for details. For their cohomology rings, the following is known.

\begin{theorem}[Theorem 1.1, \cite{AHMMS}]
Let $X(N,I)$ be the regular 
nilpotent Hessenberg 
variety determined by a lower ideal $I$ and a regular nilpotent 
element $N \in \mathfrak{g}=Lie(G)$. Then 
$$
ST(\A_I,P_1) \simeq H^*(X(N,I)).
$$
In particular, 
$$
\mbox{Poin}(X(N,I),\sqrt{x})=\prod_{i=1}^\ell (1+x+\ldots+x^{d_i^I}).
$$
\label{Hessenberg}
\end{theorem}

In \cite{AHMMS}, there are no terminology ``Solomon-Terao algebras''. Here we state the main result in \cite{AHMMS} 
in terms of the Solomon-Terao algebra. 
Theorem \ref{Hessenberg} shows that the Solomon-Terao algebra $ST(\A_I,P_1)$ is realized 
as the cohomology ring of the variety $X(N,I)$, which reminds us of that 
the Orlik-Solomon algebra is isomorphic to the cohomology ring of $M(\A_I)$. 
Note that 
$X(N,\Phi^+)=G/B$. Hence we can say that 
\textbf{the Solomon-Terao algebra generalizes the coinvariant algebra of the 
Weyl groups in the setting of hyperplane arrangements}. 

\begin{rem}
From now on, when $\K=\CC,\ \A \subset \A_W$ and $P_1$ is the same as in Theorem 
\ref{Hessenberg}.
let $ST(\A):=ST(\A,P_1)$ and $\mathfrak{a}(\A):=\mathfrak{a}(\A,P_1)$. 
It is clear that 
$P_1 \in U_2$ for any $\A$.
\end{rem}


The organization of this article is as follows. In \S2 we recall several results 
on arrangements, mainly from \cite{ST}. In \S3, we prove Theorem \ref{factorization}. 
In \S4 we investigate the Solomon-Terao algegbra for the inversion arrangements, and 
the relation to the Schubert varieties.
In 
\S5 we pose several questions related to the Solomon-Terao algebras.
\medskip

\noindent
\textbf{Acknowledgements}. The authors are 
partially supported by JSPS KAKENHI 
Grant-in-Aid for Scientific Research (B) 16H03924.

\section{Preliminaries}

In this section we collect several definitions and results, mainly from 
\cite{OT} and \cite{ST}. First, let us recall definitions on combinatorics of 
arrangements.

\begin{define}
Define the \textbf{intersection lattice} $L(\A)$ of $\A$ by 
$$
L(\A):=\{ \cap_{H \in \B} H \mid \B \subset \A\}.
$$ 
The \textbf{M\"{o}bius function} $\mu$ on $L(\A)$ is defined by, 
$\mu(V)=1$, and by $\mu(X):=-\sum_{V \supset Y \supsetneq X} \mu(Y)$. 
Then define the \textbf{Poincar\'{e} polynomial} of $\A$ by 
$$
\pi(\A;t):=\sum_{X \in L(\A)} \mu(X) (-t)^{\codim X},
$$
and define the \textbf{characteristic polynomial} of $\A$ by 
$$
\chi(\A;t):=\sum_{X \in L(\A)} \mu(X) t^{\dim X}.
$$
\label{OTbasic}	
\end{define}

By \cite{OS}, $\pi(\A;t)$ coincides with the topological Poincar\'e polynomial of 
$M(\A):=V \setminus \cup_{H \in \A}H$ when $\K=\CC$. Moreover, the presentation 
of $H^*(M(\A),\Z)$ has a presentationcan depending only on $L(\A)$, see \cite{OS} for details. 

Now let us recall several properties and results on $D(\A)$. For 
$\theta \in \Der S$, we say that $\theta$ is homogeneous of degree $d$ if 
$\deg \theta(\alpha)=d$ for all $\alpha \in V^*$ with $\theta(\alpha) \neq 0$. 
Also, let us introduce a criterion for freeness of $\A$. 

\begin{theorem}[Saito's criterion, \cite{S}]
Let $\theta_1,\ldots,\theta_\ell \in D(\A)$ be homogeneous elements. Then they form a basis for $D(\A)$ if and only if 
(1) they are $S$-linearly independent, and (2) $\sum_{i=1}^\ell \deg \theta_i=|\A|$.
\label{Saito}
\end{theorem}

The following is the most important consequence of the freeness.

\begin{theorem}[Terao's factorization, \cite{T2}]
Let $\A$ be free with $\exp(\A)=(d_1,\ldots,d_\ell)$. Then 
$$
\pi(\A;t)=\prod_{i=1}^\ell (1+d_i t).
$$
\label{Tfac}
\end{theorem}

The following plays a key role in the proof of Theorem \ref{factorization}. 

\begin{prop}[e.g., \cite{OT}, Proposition 4.12]
Let $\theta_1,\ldots,\theta_\ell \in D(\A)$. Then 
$\det (\theta_i(x_j)) \in Q(\A)S$.
\label{div}
\end{prop}

Let us introduce two sufficient conditions to check the freeness of $\A$ for 
our purpose. 

\begin{theorem}[Terao's addition-deletion theorem, \cite{T}]
Let $H \in \A,\ \A':=\A 
\setminus \{H\}$ and $\A'':=\A^H=\{H \cap L \mid L \in \A'\}$. Then two of the 
following three imply the third:
\begin{itemize}
\item[(1)]
$\A$ is free with $\exp(\A)=(d_1,\ldots,d_{\ell-1},d_\ell)$.
\item[(2)]
$\A'$ is free with $\exp(\A')=(d_1,\ldots,d_{\ell-1},d_\ell-1)$.
\item[(3)]
$\A''$ is free with $\exp(\A'')=(d_1,\ldots,d_{\ell-1})$.
\end{itemize}
\label{additiondeletion}
\end{theorem}

\begin{theorem}[Division theorem, \cite{A}, Theorem 1.1]
Let $H \in \A$. If $\A^H$ is free and $\pi(\A^H;t) \mid \pi(\A;t)$, 
then $\A$ is free.
\label{divisionthm}
\end{theorem}

Our results in this article rely on those in \cite{ST}. 
To prove Theorem \ref{polynomial}, Solomon and Terao introduced 
the Solomon-Terao complex as in Theorem \ref{ST}. 
At the same time, the structure of their complex in itself deserves our attention. 
We summarize some definitions and results from \cite{ST} below.

\begin{define}[\cite{ST}, Definition 4.5]
Let $d$ be a non-negative integer and $\A$ an arrangement in $V$.  For $X \in L(\A)$ let 
$S^X$ be the coordinate ring of $X$. We say that $h \in S^X$ is 
\textbf{non-degenerate on $X$} if 
the zero-locus of all polynomials in $\mbox{Jac}(h):=\{
\theta(h) \in S \mid \theta \in \Der S^X\}$ is contained in the 
origin of $X$. Define 
$$
U_d^X(\A):=\{f \in S_d \mid f|_X\ \mbox{ is non-degenerate on } X\}.
$$
\label{nondegenerate}
\end{define}

\begin{prop}[\cite{ST}, Section four and Corollary 3.6]
Let $\A$ be an arrangement in $V=\K^\ell$. Then 

(1)\,\,
for each $d > 0$, 
the open set 
$$
U_d(\A):=\cap_{X \in L(\A)} U_d^X(\A) \subset S_d$$ 
is non-empty, and 
the Solomon-Terao complex has a finite dimensional cohomology group
$H^i(D^*(\A), \eta)\ (i=0,\ldots,\ell)$ for all $ \eta \in U_d(\A)$.

(2)\,\,
If $\A$ is free with $\exp(\A)=(d_1,\ldots,d_\ell)$, then 
$$
\Psi(\A;x,t)=\prod_{i=1}^\ell 
(t(1+x+\cdots+x^{d_{i-1}})+x^{d_i}).
$$
\label{STbasic}
\end{prop}

By definition, the following is clear.


\begin{lemma}
Let $\eta \in U_d(\A)$. Then $\partial_{x_1}(\eta),\ldots,
\partial_{x_\ell}(\eta)$ is an $S$-regular sequence.
\label{regular}
\end{lemma}

\noindent
\textbf{Proof}. 
Apply the definition of $U_d^X(\A)$ when $X=V$. \owari
\medskip

In arrangement theory, for $H \in \A$, we often consider the deletion 
$\A':=\A \setminus \{H\}$ and the restriction $\A'':=\A^H:=\{L \cap H \mid 
L \in \A'\}$ together to obtain 
the information of $\A$. However, it is not easy to see whether 
$\eta \in U_d(\A)$ is also contained in $U_d(\A')$ or not. For the restriction, we have the 
following.

\begin{lemma}
Let $\eta \in U_d(\A)$ and $H \in \A$. Then 
$\eta \in U(\A \setminus \{H\}))$, and $\eta|_H \in U_d(\A^H)$.
\label{restdel}
\end{lemma}

\noindent
\textbf{Proof}. 
If $X \in L(\A \setminus \{H\})$, then $X \in L(\A)$ since $L(\A \setminus \{H\}) \subset L(\A)$. 
If $X \in L(\A^H)$, then so is $X \in L(\A)$ since $L(\A^H) \subset L(\A)$. \owari
\medskip

Hence for $\eta \in U(\A)$, we have the following two maps:
\begin{eqnarray}
&F_1&:ST(\A \setminus \{H\},\eta) \stackrel{\cdot \alpha_H}{\rightarrow} ST(\A,\eta),\label{eq100}\\
&F_2&:ST(\A,\eta) \rightarrow ST(\A^H,\eta|_H).\label{eq101}
\end{eqnarray}
Also, it is clear that $F_2 \circ F_1=0$ and $F_2$ is surjective.
To investigate a general property of $\Psi(\A;x,t)$, we use the following.

\begin{prop}[\cite{ST}, Proposition 4.4]
For $H \in \A$, consider the boundary map 
$\partial_p^H:D^p(\A) \rightarrow D^{p-1}(\A)\ (p=1,\ldots,\ell)$ defined by 
$$
\partial_p^H(\theta)(f_2,\ldots,f_p):=
\displaystyle \frac{\theta(\alpha_H,f_2,\ldots,f_p)}{\alpha_H}
$$
for $\theta \in D^p(\A)$. Then the complex $(D^*(\A),\partial_p^H)$ is exact. In particular, 
$$
\sum_{p=0}^\ell \mbox{Hilb}(D^p(\A);x)(-x)^{\ell-p}=0.
$$
\label{acyclic}
\end{prop}

\begin{prop}

Assume that $\A$ is tame, i.e., 
$$
\mbox{pd}_S D(\A)^p \le \ell- p\ (p=0,\ldots,\ell).
$$
Then 
$H^i(D^*(\A), \eta)=0$ for $i \neq 0$. 

\label{tame}
\end{prop}

\noindent
\textbf{Proof}. Theorem 5.8 in \cite{OT2} states that for the complex
$$
0 \rightarrow C^0 \rightarrow \cdots \rightarrow C^\ell \rightarrow 0
$$
of $S:=\K[x_1,\ldots,x_\ell]$-modules, the cohomology group 
$H^k$ of this complex vanishes if 
$$
\mbox{pd}_S C^p < \ell+p-k\ (\forall p).
$$
Now $C^p=D^{\ell-p}(\A)$. Hence 
$$
\mbox{pd}_S C^p \le p= \ell+(p-\ell)<\ell+(p-k)
$$
for all $p$ and $k \neq \ell$. Hence 
$H^k(D^*(\A),\partial_*)=0$ for $k \neq \ell$.\owari
\medskip

The following result in commutative ring theory play the key roles in the proof of Theorem \ref{factorization}. 

\begin{theorem}[e.g., \cite{Sm}, Theorem 6.5.1]
Let $S=\K[x_1,\ldots,x_\ell]$ and let $t_1,\ldots,t_\ell$ be an $S$-regular sequence.
Let $I=(t_1,\ldots,t_\ell)$ and assume that 
$t_i=\sum_{j=1}^\ell a_{ij} x_j$ for $i=1,\ldots,\ell$. Then $\Delta:=\det (a_{ij})$ is a 
$\K$-basis for $\soc(S/I)$.
\label{Smith}
\end{theorem} 


\begin{theorem}[e.g., \cite{Sm}, Theorem 6.7.6]
Let $t_1,\ldots,t_\ell$ be homogeneous polynomials with $\deg t_i=d_i\ (i=1,\ldots,\ell)$, and let $R=S/I$ for an ideal $I =(t_1,\ldots,t_\ell)$. 
If $\dim_\K R< \infty$, then $t_1,\ldots,t_\ell$ is an $S$-regular sequence, and 
$$
\mbox{Hilb}(R;x)=\prod_{i=1}^\ell (1+x+\cdots+x^{d_i-1}).
$$
\label{Smith3}
\end{theorem}

\section{Proof of Theorems \ref{poin2} and \ref{factorization}}

From now on let us fix 
$\eta \in U_d(\A)$ unless otherwise specified. 
First we prove Theorem \ref{poin2}. 
For that, let us show the following proposition essentially proved in \cite{ST}.


\begin{prop}
For $\eta \in U_2(\A)$ and an arbitrary arrangement $\A$, 
$$
\sum_{p=0}^{\ell} (-x)^p \mbox{Hilb}(H^p(D^*(\A),\eta);x)=\Psi(\A;x,1).
$$
\label{poincplx}
\end{prop}

\noindent
\textbf{Proof}. 
By the results in \S2, we know that 
\begin{eqnarray*}
\Psi(\A;x,1)&=&\sum_{p=0}^\ell 
\mbox{Hilb}(D^p(\A);x)((1-x)-1)^p\\
&
=&
\sum_{p=0}^\ell \mbox{Hilb}(D^p(\A);x)(-x)^p\\
&=&
\sum_{p=0}^\ell \mbox{Hilb}(H^p(D^*(\A));x)(-x)^p.
\end{eqnarray*}
Here we used the fact that 
$\partial_p$ is of degree one since $\eta \in U_2(\A)$, and the finite dimensionality of 
$H^p(D^*(\A))$ by Theorem \ref{ST}. \owari
\medskip



\noindent
\textbf{Proof of Theorem \ref{poin2}}. 
Combine Proposition \ref{poincplx} with Proposition \ref{tame} 
and 
the properties of $\pi(\A;t)$.\owari
\medskip

Next let us prove Theorem \ref{factorization}. 
\medskip

\noindent
\textbf{Proof of Theorem \ref{factorization}}. 
(1)\,\,
Let us show 
that $ST(\A, \eta)$ is a complete intersection. Let $\theta_1,\ldots,\theta_\ell$ be 
a homogeneous basis for $D(\A)$, and let $\theta_i(\eta)=:f_i$. By Theorem \ref{ST}, we know that $ST(\A, \eta)=S/\mathfrak{a}(\A, \eta)$ is a finite dimensional 
$\K$-algebra, and $\mathfrak{a}(\A, \eta)$ is generated by $\ell$-homogeneous 
polynomials $f_1,\ldots,f_\ell$. Then $f_1,\ldots,f_\ell$ form 
a regular sequence by Theorem \ref{Smith3}, and hence $ST(\A, \eta)$ is a complete intersection. 
On the Hilbert series, apply Theorem \ref{Smith3}. 

(2)\,\,
To prove (2), let us prove the following proposition.

\begin{prop}
Let $\eta \in U_d(\mathcal A)$ and let $\eta_{ij}=\partial_{x_i}\partial_{x_j} \eta$ for $1 \leq i,j \leq \ell$.
The element $Q(\A) \det(\eta_{ij})$ is contained in $\soc(ST(\A,\eta))$.
Moreover, if $ST(\A,\eta)$ is a complete intersection, then $Q(\A) \det(\eta_{ij})$ is a $\mathbb K$-basis for $\soc({ST(\A,\eta)})$. 
\label{socleelement}
\end{prop}

\noindent
\textbf{Proof}. 
By Lemma \ref{regular}, $\eta \in U_d(\A)$ implies that 
$\partial_{x_1}(\eta),\ldots,\partial_{x_\ell}(\eta)$ form an $S$- 
regular sequence. Let $\eta_i:=\partial_{x_i}(\eta)$. Since $\eta_i$ is homogeneous, 
it holds that 
$$
\eta_i=\sum_{j=1}^\ell \eta_{ij} x_j
$$
up to non-zero scalar. 
Hence Theorem \ref{Smith} implies that 
$\zeta:=\det(\eta_{ij})$ is a $\K$-basis for $\soc(ST(\emptyset,\eta))$. 
Since $\emptyset \subset \A$, we have an $S$-morphism 
$$
F:ST(\emptyset,\eta) \rightarrow ST(\A,\eta)$$
defined by multiplying $Q(\A)$, which is well-defined by definition of 
the Solomon-Terao algebra.
%
Since 
$$
xQ(\A) \zeta=x F(\zeta)=F(x\zeta)=F(0)=0
$$
for any $x \in S_+$, it holds that 
$Q(\A) \zeta \in \soc(ST(\A,\eta))$.

Now, we assume that $ST(\A,\eta)$ is a complete intersection,
and let us show that $Q(\A) \zeta \neq 0$ 
in $ST(\A,\eta)$. 
Let $\theta_1,\ldots,\theta_\ell \in D(\A)$ such that 
$f_1:=\theta_1(\eta),\ldots,f_\ell:=\theta_\ell(\eta)$ is an $S$-regular sequence belonging to 
$\mathfrak{a}(\A,\eta)$. Let 
$\theta_i=\sum_{j=1}^\ell f_{ij} \partial_{x_j}$. 
Then we have
$$
f_i=\theta_i(\eta)=
\sum_{k=1}^\ell f_{ik} \eta_k=
\sum_{k=1}^\ell (\sum_{j=1}^\ell f_{ik} \eta_{kj} x_j)
=
\sum_{j=1}^\ell (\sum_{k=1}^\ell f_{ik} \eta_{kj}) x_j
$$
up to nonzero scalor. Hence Theorem \ref{Smith} implies that 
$$
\det (\sum_{k=1}^\ell f_{ik}\eta_{kj})=\det (f_{ij}) \det (\eta_{ij})=\det (f_{ij}) \zeta
$$
is a $\K$-basis for $\soc(ST(\A,\eta))$. 
%
%
In particular,  
$\det (f_{ij}) \zeta$ is not zero in $ST(\A,\eta)$. By Proposition \ref{div}, 
we know that $\det (f_{ij})$ is divisible by $Q(\A)$. 
In other words, there is $g \in S$ such that 
$\det (f_{ij}) \zeta=g Q(\A) \zeta$ in $ST(\A,\eta)$. Since both are 
elements of $\soc(ST(\A,\eta))$, $g$ is a nonzero-scalor. \owari
\medskip

\noindent
\textbf{Proof of Theorem \ref{factorization} (2), continued}. 
Now assume that $\eta \in U_d(\A)$ and $ST(\A,\eta)$ is a complete 
intersection ring. Let $\theta_1,\ldots,\theta_\ell \in D(\A)$ be derivations such that $\theta_1(\eta),\ldots,\theta_\ell(\eta)$ form an $S$-regular sequence. 
Let $\zeta:=\det (\eta_{ij})$. 
We show that $\theta_1,\ldots,\theta_\ell$ form a basis for 
$D(\A)$ by Saito's criterion. For $\theta_i =\sum_{j=1}^\ell f_{ij} \partial_{x_j}\in D(\A)$, 
by Proposition \ref{socleelement}, we know that 
$Q(\A) \zeta=\det (f_{ij}) \zeta$ is a $\K$-basis for $\soc(ST(\A,\eta))$. Thus 
$\theta_1,\ldots,\theta_\ell$ are $S$-independent. Moreover, 
$\deg \det (f_{ij})= \sum_{i=1}^\ell \deg \theta_i=\deg Q(\A)=|\A|$. Hence
%
Saito's criterion 
implies that $\A$ is free. On the Hilbert series, apply Theorem \ref{Smith3}. \owari
\medskip


%

In the following low-dimensional cases, we can always apply some of results above.

\begin{prop}[Two and three dimensional case]
(1)\,\,
Assume that $\ell=2$. Then all non-empty arrangements are free and hence tame. In particular, 
$$
\mbox{Hilb}(ST(\A, \eta);x)=(1+x+\cdots+x^{d-1})(1+x+\cdots+x^{d+|\A|-3})
$$
for $\eta \in U_d(\A)$.

(2)\,\,
Assume that $\ell=3$. Then all arrangements are tame.
Hence $
\mbox{Hilb}(ST(\A,\eta);x)=\Psi(\A;x,1)$ 
for $\eta \in U_2(\A)$.
\label{23}
\end{prop}

\noindent
\textbf{Proof}. 
(1)\,\,
It is famous that $\A$ is free with $\exp(\A)=(1,|\A|-1)$ 
when $\ell=2$. Hence 
Theorem \ref{factorization} completes the proof.

(2)\,\, Since $D^p(\A)$ is reflexive, their projective dimension is 
at most $1$. Hence the only case we have to check the tameness 
is whether $\mbox{pd}_S D^3(\A) \le 0$. This is true since $D^3(\A) \simeq S$. For the rest, apply Theorem \ref{poin2}. \owari
\medskip


In general, it is not easy to compute Solomon-Terao polynomials and 
$\mbox{Hilb}(ST(\A,\eta);x)$ when $\A$ is not free.

\begin{example}
Let $\A$ be defined as $xyz(x+y+z)=0$. Then we may compute 
$\mbox{Hilb}(ST(\A);x)=1+3x+5x^2+4x^3+x^4=(1+x)(1+2x+3x^2+x^3)$.
\label{4}
\end{example}



It is known that $\A$ is not necessarily free even if $\pi(\A;t)=\prod_{i=1}^\ell 
(1+d_i t)$ with $d_1,\ldots,d_\ell \in \Z$. For the Solomon-Terao algebras, we do not know 
any such examples. Based on several computations, we pose the following conjectures.

\begin{conj}
(1)\,\,
$\A$ is free if and only if $$
\mbox{Hilb}(ST(\A, \eta);x)=\prod_{i=1}^\ell (1+x+\cdots+x^{d_i})
$$
for some integers $d_1,\ldots,d_\ell$.

(2)\,\,
$\A$ is free if and only if $
\mbox{Hilb}(ST(\A, \eta);x)$ is palindromic, i.e., for 
$
\mbox{Hilb}(ST(\A, \eta);x)=\sum_{i=0}^n a_i x^i$ with $a_n \neq 0$, it holds that 
$a_i=a_{n-i}$ for all $i$.

\label{conj3}
\end{conj}

The ``only if'' part of Conjecture \ref{conj3} is surely true by 
Theorem \ref{factorization}. Let us check Conjecture \ref{conj3} when $\A$ is not free but 
$\pi(\A;t)$ splits over $\Z$ for the following case.

\begin{example}
Let $$
\A:=\{x(x^2-y^2)(x^2-2y^2)(y-z)z=0\}.
$$
It is easy to check that 
$\pi(\A;t)=(1+t)(1+3t)^2,
$ but $\A$ is not free (hence the factorization of $\pi(\A;t)$ is not the 
sufficient condition for the freeness). In this case, let us compute 
$\mbox{Hilb}(ST(\A,\eta);x)$.
Let $\eta:=x^2+y^2+z^2$. Then we can compute that 
$$
\mathfrak{a}(\A,\eta)=(x^2+y^2+z^2,\ z^3-yz^2,\ 
y^6-y^5z,\ y^6+3y^4z^2).
$$
Hence
\begin{eqnarray*}
\mbox{Hilb}(ST(\A,\eta);x)&=&
1+3x+5x^2+6x^3+6x^4+6x^5+4x^6+x^7\\
&=&(1+x)^2(1+x+2x^2+x^3+2x^4+x^5),
\end{eqnarray*}
which does not split into the form $\prod (1+x+\cdots+x^{d_i-1})$. Thus this example does not give a counter example to Conjecture \ref{conj3}. 
\label{notsplit}
\end{example}

We give one result related to Conjecture \ref{conj3} as follows.

\begin{prop}
Let $\A$ be an arrangement in $\K^3$ and $\eta \in U_2(\A)$. Assume that 
$$
\mbox{Hilb}(ST(\A,\eta);x)=\prod_{i=1}^3 (1+x+\cdots+x^{d_i})
$$
for $d_1=1$ and some $d_2,d_3 \in \Z$. If there is $H \in \A$ such that 
$$
\mbox{Hilb}(ST(\A^H,\eta|_H);x)=\prod_{i=1}^2 (1+x+\cdots+x^{d_i}),
$$
then $ST(\A,\eta)$ is a complete intersection.
\label{point}
\end{prop}

\noindent
\textbf{Proof}. 
By Lemma \ref{restdel}, $\eta|_H \in U_2(\A^H)$. 
By Proposition \ref{23} (1), $ST(\A^H,\eta|_H)$ is always a complete intersection, and 
$\A$ is tame by Proposition \ref{23} (2). Let us compute the coefficients of $\pi(\A;t)
=1+|\A|t+b_2t^2+b_3t^3$ in terms 
of $\mbox{Hilb}(ST(\A,\eta);x)$. By Theorem \ref{poin2}, 
$$
\Psi(\A;1,1)=\pi(\A;1)=
\mbox{Hilb}(ST(\A,\eta);1)=1+|\A|+b_2+b_3
=2(d_2+1)(d_3+1).
$$
Also, since $\A$ is central, $\pi(\A;t)$ is divisivle by $1+t$, see \cite{OT} for example. Hence 
$$
\pi(\A;-1)=1-|\A|+b_2-b_3=0.
$$
By these two we can compute $b_2=d_2+d_3+d_2d_3$. Again by Proposition \ref{point} (2), 
$\pi(\A^H;t)=(1+t)(1+d_2t)$. By Theorem \ref{divisionthm} when $\ell=3$,  
$b_2=d_2+|\A^H|(|\A|-|\A^H|)=d_2+(d_2+1)d_3$. 
Hence $\A$ is free, and Theorem \ref{factorization} implies that $ST(\A,\eta)$ is 
a complete intersection. \owari
\medskip

\section{Inversion arrangements and 
Schubert varieties}

In this section we use the notation in \S1.2, i.e., $\K=\CC$ and 
$\eta=P_1$. 

\begin{define}
Let $\Phi$ be the root system with respect to the Weyl group $W$ and fix a 
positive system $\Phi^+$. 
For $w \in W$, define $\A_w:=\{\alpha=0\mid \alpha \in \Phi^+,\ 
w \alpha \ \mbox{is a negative root}\}$, 
which is called the \textbf{inversion arrangement}. 
\label{inv}
\end{define}

Also for $w \in W$, we have the Schubert variety 
$Y_w:
=\overline{BwB}$. 
Now let us check the freeness of inversion arrangements. 
For details of them, see \cite{Slof}. 

\begin{theorem}[Theorem 3.3, \cite{Slof}]
Let $Y_w$ be the Schubert variety determined by $w \in W$. Then 
$Y_w$ is rationally smooth if and only if 
$\A_w$ is free with $\exp(\A_w)=(d_1^w,\ldots,d_\ell^w)$, and 
$\prod_{i=1}^\ell (1+d_i^w)=|[e,w]|$, the number of elements between $e$ and $w$ in 
the Bruhat order. Moreover, 
$$
\mbox{Poin}(Y_w;\sqrt{x})=\prod_{i=1}^\ell (1+x+\ldots+x^{d_i^w}).
$$
\label{inversion}
\end{theorem}

Hence by Theorem \ref{factorization}, we have the following:

\begin{cor}

In the setup of Theorem \ref{inversion}, 
$$
\mbox{Poin}(Y_w;\sqrt{x})=\mbox{Hilb}(ST(\A_w,P_1);x).
$$

\label{poincare}

\end{cor}

The claim above suggests a correspondence between inversion arrangements 
and Schubert varieties. We have the Solomon-Terao algebra on one side 
and the cohomology ring on the other. However, they are not isomorphic 
as algebras in general. 

\begin{prop}
Let $w=(4123) \in S_4$. Then $\A_w$ is defined by 
$$
(x_1-x_2)(x_1-x_3)(x_1-x_4)=0
$$
in $\CC^4$. Then 
$$
H^*(Y_w,\CC) \not \simeq ST(\A_w,P_1)
$$
as rings.
\label{Yw}
\end{prop}

\noindent
\textbf{Proof}. By the computation of Schubert polynomials, it holds that 
$$
H^*(Y_w,\CC) \simeq \CC[x_1,x_2,x_3,x_4]/(f_1,f_2,f_3,f_4),
$$
where 
\begin{eqnarray*}
f_1&=&x_1+x_2+x_3+x_4,\\
f_2&=&(x_1+x_2+x_3)^2,\\
f_3&=&x_2x_3++x_1x_3,\\
f_4&=&x_1x_2.
\end{eqnarray*}
Hence $H^*(Y_w,\CC)$ has a non-zero element $x_1+x_2+x_3$ of degree one 
such that $(x_1+x_2+x_3)^2=0$. 
On the other hand, we can check by the direct computation that there are no such elements of degree one in $ST(\A_w,P_1)$. Hence they are not 
isomorphic.\owari
\medskip

So 
the statement of Theorem \ref{Hessenberg} for Hesssenberg varieties does not hold 
for
Schubert varieties in general, though we have Corollary \ref{poincare}. However, since 
the Solomon-Terao algebra depends on the choice of $\eta \in U_2(\A)$, we may pose the following problem.

\begin{problem}
Are there any $\eta \in U_2(\A_w)$ such that 
$$
H^*(Y_w) \simeq ST(\A_w,\eta)
$$
as rings?
\label{conj}
\end{problem}

By Theorem \ref{Hessenberg}, Problem \ref{conj} has a positive answer for $\eta=P_1$ when 
$w$ is the longest element in $W$, and 
for some special $w$ by Theorem \ref{Hessenberg}. 



\section{Questions, problems and conjectures}

In this section we collect several problems and conjectures related to Solomon-Terao algebras. 
We assume that $\eta \in U_2(\A)$ unless otherwise specified. 
The most important problem is the following.

\begin{question}
Are there any topological meaning of $\Psi(\A;x,t)$ and 
$\mbox{Hilb}(ST(\A,\eta);x)$? Also, are there any 
common property among them?
\label{important}
\end{question}

For geometric meaning, by Theorem \ref{Hessenberg}, 
we can say that $\Psi(\A;x,1)=\mbox{Hilb}(ST(\A,P_1),x)$ is the Poincar\`{e} 
polynomial of the regular nilpotent Hessenberg varieties when $\A$ is the ideal 
arrangement. For general properties, 
we can say the following, which is essentially proved in Proposition 5.4, \cite{ST}. 

\begin{prop}
If $\A \neq \emptyset$, then $x+t$ divides $\Psi(\A;x,t)$. Therefore, 
$1+x$ divides $\mbox{Hilb}(ST(\A,\eta);x)$ when 
$\A$ is tame. 
\label{xt}
\end{prop}

\noindent
\textbf{Proof}. By Proposition 
\ref{acyclic}, 
$$
\sum_{p=0}^\ell \mbox{Hilb}(D^p(\A),x)(-x)^{\ell-p}=0.
$$
On the other hand, 
$$
\sum_{p=0}^\ell \mbox{Hilb}(D^p(\A),x)(-x)^{\ell-p}=\Psi(\A;x,-x)=0.
$$
Since $\Psi(\A;x,-x)$ is a polynomial in $x$, we complete the proof. \owari
\medskip

Another question is to ask whether 
we can compute $\mbox{Hilb}(ST(\A,\eta);x)$ inducively as for $\pi(\A;t)$ or not. 
For $\pi(\A;t)$, let $H \in \A,\ \A':=\A 
\setminus \{H\}$ and $\A'':=\A^H$. Then it holds that 
$$
\pi(\A;t)=\pi(\A';t)+t\pi(\A'';t),
$$
which is called the deletion-restriction formula. 

\begin{question}
Is there a deletion-restriction type formula for $\mbox{Hilb}(ST(\A,\eta);x)$? i.e., 
for the triple $(\A,\A',\A'')$, is there a formula between
$
\mbox{Hilb}(ST(\A,\eta);x),\ 
\mbox{Hilb}(ST(\A',\eta);x)$, and $\mbox{Hilb}(ST(\A^H,\overline{\eta});x)$?
Also, what about the same question for the Solomon-Terao polynomials 
$\Psi(\A;x,t)$?
\label{DR}
\end{question}

The naive generalization of the deletion-restriction formula does not 
work well in general as follows.

\begin{example}
The most simply-minded generalization is as follows:
\begin{equation}
\mbox{Hilb}(ST(\A,\eta);x)=
\mbox{Hilb}(ST(\A',\eta);x)+x^{|\A|-|\A^H|}\mbox{Hilb}(ST(\A^H,\eta|_H);x).
\label{eq1}
\end{equation}
This does not hold true in general. 
Let $\A$ be an arrangement defined by
$xyz(x+y+z)=0$. Then we have  
$$
\Psi(\A;x,1)=1+3x+5x^2+4x^3+x^4.
$$
Let $H=\{x+y+z=0\} \in \A$. Then it is easy to check that $\A':=\A \setminus \{H\}$ and $\A^H$ are both free with 
exponents $(1,1,1)$ and $(1,1)$. Thus by Theorem \ref{factorization}, it holds that 
\begin{eqnarray*}
\Psi(\A';x,1)&=&(1+x)^3,\\
\Psi(\A^H;x,1)&=&(1+x)(1+x+x^2).
\end{eqnarray*}
Hence the deletion-restriction does not hold for $\Psi(\A;x,1)$, neither 
does $\Psi(\A;x,t)$.

By Theorems \ref{factorization}, \ref{additiondeletion} and 
\ref{divisionthm}, the formula (\ref{eq1}) holds true 
either when (a) $\A$ and $\A'$ are free, or 
(b) $\A^H$ is free and $\pi(\A^H;t)$ divides $\pi(\A;t)$.
\end{example}

\begin{question}
Are there any similarity between the Solomon-Terao polynomials or 
$\mbox{Hilb}(ST(\A,\eta);x)$, and 
Poincar\'{e} polynomials of arrangements for the triple $(\A,\A',\A^H)$? For example, 
the division 
$$
\pi(\A^H;t) \mid \pi(\A;t)$$
implies 
$$
\mbox{Hilb}(ST(\A^H,\eta|_H);x) \mid \mbox{Hilb}(ST(\A,\eta);x)
$$
and vice versa?
\label{division}
\end{question}

Again Problem \ref{division} is true for free cases.

\begin{prop}
Assume that, either (a) $\A$ and $\A'$ are free, or 
(b) $\A^H$ is free and $\pi(\A^H;t)$ divides $\pi(\A;t)$. Then 
$\mbox{Hilb}(ST(\A^H);x) \mid \mbox{Hilb}(ST(\A);x)$. 
\label{DRtrue}
\end{prop}

\noindent
\textbf{Proof}. 
Apply Theorems \ref{factorization}, \ref{additiondeletion} and 
\ref{divisionthm}. \owari
\medskip

By the definition of $ST(\A)$, we can treat some 
finite-dimensional algebras associated with tame arrangements. Hence it will be 
meaningful to give names to arrangements depending on the properties 
as follows:

\begin{define}
Let $\A$ be tame. Then we say that 

(1)\,\,
$(\A,\eta)$ is \textbf{of complete intersection} if $ST(\A,\eta)$ is a complete intersection ring, 

(2)\,\,
$(\A,\eta)$ is \textbf{Gorenstein} if $ST(\A,\eta)$ is a Gorenstein ring, or equivalently, 
$ST(\A,\eta)$ is a Poincar\'{e} duality algebra,

(3)\,\,
$(\A,\eta)$ is \textbf{SLP} if $ST(\A,\eta)$ has a strong Lefschetz element (see 
\cite{lef} for details), and 



(4)\,\,
$\A$ is \textbf{ST-finite} if $H^i(D^*(\A),\eta)=0$ for $i \neq 0$.
\label{names}
\end{define}

\begin{rem}
Note that the properties in Definition \ref{names} depend on the choice of 
$\eta \in U_2(\A)$ in general. 
See the following example. So to investigate which properties depends only on $\A$ may be 
an interesting question. 
For example, when $\eta=Q(\A)$, it follows by definition of $D(\A)$ that 
$\mathfrak{a}(\A,Q(\A))=SQ(\A)$ for a non-empty $\A$. Hence ST-finiteness depends on the choice of $\eta \in S_d$.
\label{rem1}
\end{rem}

\begin{example}
Let $\A:=\{xy(x+y)=0\}$. Then 
$D(\A)$ has a basis $\theta_E,y(x+y)\partial_y$.  Let 
$\eta_0:=x^4+y^4 \in U_4(\A)$ and $\eta:=\sum_{i=0}^4 a_i x^iy^{4-i}$ for 
$a_i \in \R$. We show that for generic $a_0,\ldots,a_4$, two Solomon-Terao algebras 
$ST_0:=ST(\A,\eta_0)$ and $ST:=ST(\A,\eta)$ are not isomorphic as 
graded rings. Note that 
$\mathfrak{a}(\A,\eta_0)_{\le 1}=\mathfrak{a}(\A,\eta)_{\le 1}=0$. Hence if a 
graded ring isomorphism $\psi:ST_0 \rightarrow ST$ exists, then it is induced from 
a graded ring isomorphism $\varphi:S \rightarrow S$ since they are generated by degree one part. Let $\varphi(x)=\alpha x+\beta y,\ 
\varphi(y)=\gamma x+\delta y$. Since $\psi$ is induced from $\varphi$, it holds that 
$\varphi(\eta_0) \in (\eta)$. We can see that this can be expressed as a closed condition. 
Hence generically, there are no $\varphi$ which induces $\psi$. Hence $ST_0 \not \simeq 
ST$.
\label{notisom}
\end{example}



By Theorem \ref{factorization}, a complete intersection property 
is same as the freeness independent of the choice of $\eta \in U_d(\A)$. 
Now we may pose the 
following natural problems:

\begin{question}
Give a sufficient, or equivalence condition for the Gorenstein or ST-finite arrangement. 
\label{namesproblem}
\end{question}

Also, not all arrangements are Gorenstein.

\begin{prop}
The arrangement $xyz(x+y+z)=0$ is not Gorenstein.
\label{notGorenstein}
\end{prop}

\noindent
\textbf{Proof}. 
As seen in Example \ref{4}, 
$
\mbox{Hilb}(ST(\A);x)=1+3x+5x^2+4x^3+x^4$, which is not palindromic. 
Hence this arrangement is not Gorenstein.
\owari
\medskip


To ask 
the top degree of the nonzero part of 
the Solomon-Terao algebra is a natural question. 
By some computation, we conjecture the following: 

\begin{conj}
Let $\eta \in U_d(\A)$ and define $r:=\max\{n \mid 
ST(\A,\eta)_n \neq (0)\}$. Then 

(1)\,\, $r=|\A|+\ell(d-2)$, and 

(2)\,\,
$\dim_\K ST(\A,\eta)_r=1$.

\label{socdegA}
\end{conj} 

Under some generic condition, we can give a partial affirmative answer to Conjecture \ref{socdegA}.

\begin{theorem}
Let $\mathcal A$ be an arrangement of linear hyperplanes.
For a generic $\eta \in U_d(\mathcal A)$, the element $Q(\mathcal A) \det (\eta_{ij})$ is a nonzero element of $\mathrm{soc}(ST(\mathcal A,\eta))$, in particular, $\dim_{\mathbb K} ST(\A,\eta)_{|\mathcal A|+\ell(d-2)} \geq 1$, where $\eta_{ij}=\partial_{x_i}\partial_{x_j} \eta$ for all $1 \leq i,j \leq \ell$. \label{socdegAgeneric}
\end{theorem}

\noindent
\textbf{Proof}. 
It is known that there exists a supersolvable arrangement $\B$ containing $\A$. See the proof of Proposition 3.5 in [AK] for example.
By the genericity, we may assume $\eta \in U_d(\A) \cap U_d (\B)$.
We already see in the Proposition \ref{socleelement} that $Q(\A) \det (\eta_{ij})$ is an element of $\soc(ST(\A,\eta))$. We claim that this is non-zero. By Theorem \ref{factorization}, the Solomon-Terao algebra of $\B$ is a complete intersection. 
Consider the $S$-morphism $F:ST(\A,\eta) \to ST(\B,\eta)$ sending $\alpha$ to $(Q(\B)/Q(\A))\alpha$,
which is well-defined since $(Q(\B)/Q(\A))D(\A) \subset D(\B)$.
Since $\B$ is free,
$F(Q(\A) \det(\eta_{ij}))=Q(\B)\det(\eta_{ij})$ is a non-zero element in $ST(\B,\eta)$ by Proposition \ref{socleelement}.
Hence $Q(\A) \det (\eta_{ij})$ is also a non-zero element in $ST(\A,\eta)$ as desired.
\owari
\medskip

\begin{problem}
Assume that $\mbox{char}(\K)=0$ and $(\A,\eta)$ is Gorenstein. 
Then by Lemma 3.74 in \cite{lef}, there is a homogeneous polynomial 
$h_\A \in \K[y_1,\ldots,y_\ell]$ such that 
$$
ST(\A) \simeq Q/\mbox{Ann}_Q(h_\A).
$$
Here $Q=\K[\partial_{y_1},\ldots,\partial_{y_\ell}]$ and 
$$
\mbox{Ann}_Q(h_\A):=\{ \epsilon \in Q \mid \epsilon(h_\A)=0\}.
$$
See \cite{lef} for details. Then determine $h_\A$. 
\label{ann}
\end{problem}

For Hessenberg varieties, 
Theorem 11.3 in \cite{AHMMS} explicitly determined 
$h_\A$. Based on Theorem \ref{socdegAgeneric} and results in Problem \ref{ann}, we can show some genericity result on 
Gorenstein arrangement.

\begin{theorem}
Let $\mbox{char}(\K)=0$, $\A$ be tame with a free arrangement $\B$ containing $\A$. 

(1)\,\,Assume that 
$ST(\A,\eta)$ is Gorenstein for $\eta \in U:=U_2(\A)$. 
Then 
there is an open set $V \subset U_2(\A)$
containing $\eta$ 
such that $ST(\A,\eta')$ is Gorenstein for all $\eta' \in V$. 

(2)\,\, 
Assume that 
$(\A,\eta)$ is SLP for $\eta \in U=U_2(\A) \cap U_2(\B)$. Then 
there is an open set $V \subset U_2(\A) \cap U_2(\B)$ containing $\eta$ 
such that $ST(\A,\eta')$ is SLP for all $\eta' \in V$. 
\label{Ggeneric}
\end{theorem}

\noindent
\textbf{Proof}. 
(1)\,\,
Let $\zeta \in U$ and denote 
$ST_\zeta:=ST(\A,\zeta)$. 
By Theorem \ref{poin2}, it holds that 
$$
\mbox{Hilb}(ST(\A,\eta);x)=
\mbox{Hilb}(ST(\A,\zeta);x).
$$
Hence for the top degree $r$ of $ST_\zeta$, it holds that 
$\dim_\K (ST_\zeta)_r=1$. Fix an isomorphism 
$[\ \ \ ]:(ST_\zeta)_r \rightarrow \K$. Define a polynomial 
$F_\zeta(y_1,\ldots,y_\ell)$ by 
$$
F_\zeta:=\displaystyle \frac{1}{r!}
\left[
\left(
\sum_{i=1}^\ell x_iy_i
\right)^r\ 
\right] \in \K[y_1,\ldots,y_\ell].
$$
Let $STG_\zeta:=Q/\mbox{Ann}_Q (F_\zeta)$ be the Gorenstein algebra. We show that there is a surjection 
$ST_\zeta \rightarrow STG_\zeta$ by sending $x_i$ to $\partial_{y_i}$. 
To show it, it is sufficient to show that 
$\mathfrak{a}_\zeta:=\mathfrak{a}(\A,\zeta) 
\subset \mbox{Ann}_Q(F_\zeta)$ regarding $S =Q$. 

Define $$
E:=\sum_{n=0}^\infty \displaystyle \frac{1}{n!}(\sum_{i=1}^\ell x_i y_i)^n.
$$
By definition, $[E]=F_\zeta$ in $ST_\zeta$. Also, note that 
$$
f(\partial_{y_1},\ldots,\partial_{y_\ell})E=
f(x_1,\ldots,x_\ell)E.
$$
Now let $f(x_1,\ldots,x_\ell) \in \mathfrak{a}_\zeta$. Then 
\begin{eqnarray*}
f(\partial_{y_1},\ldots,\partial_{y_\ell})F_\zeta&=&
f(\partial_{y_1},\ldots,\partial_{y_\ell})[E]\\
&=&
[f(\partial_{y_1},\ldots,\partial_{y_\ell})E]\\
&=&
[f(x_1,\ldots,x_\ell)E].
\end{eqnarray*}
Thus $f(\partial_{y_1},\ldots,\partial_{y_\ell})F_\zeta=0$ implies that 
$f \in \mbox{Ann}_Q(F_\zeta)$. Thus we have the surjection 
$g_\zeta:ST_\zeta \rightarrow STG_\zeta$.

Now let us show that $g_\zeta$ is injective too for generic $\zeta$. Let 
$K_\zeta$ be the kernel of $g_\zeta$. 
By the assumption, $g_\eta$ is injective, equivalently, 
$K_\eta=(0)$. Thus so is $K_\zeta=(0)$ at the neighborhood $V$ of $\eta$. 

(2)\,\,
When $\A$ is essential, it holds that $ST(\A,\eta)_1=S_1$. 
If there exists a strong Lefschetz element $\alpha$ for $ST(\A,\eta)$, then 
we can define a global morphism 
$$
\cdot \alpha:ST(\A,\zeta) \rightarrow ST(\A,\zeta)
$$
for generic $\zeta$. 
Then the same proof as (1) 
implies that the strong Lefschetz property is also generic.\owari
\medskip

Also, we have the following simple but important problem.

\begin{problem}
Consider $ST(\A,\eta)$ for non-tame $\A$. We may define 
it without the assumption of the tameness, but as far as we know, it seems difficult to 
say properties of $ST(\A,\eta)$ in that case.
\end{problem}

\end{document}